\magnification\magstep1
\hsize=16.8truecm
\vsize=25.4truecm
\baselineskip=12pt

\tolerance=1200
%
%
%
\long\def\title#1{\parindent 0pt
{\baselineskip 24pt\tolerance=10000
\vglue 3.9truecm
\noindent{#1}}\par}
\long\def\author#1{\vskip 12pt\begingroup{
\tolerance=10000\parindent=3truecm #1}\par\endgroup}

%
\countdef\sectno=11
\sectno = 0
\countdef\sbsectno=12
\sbsectno = 0
\countdef\ssbsectno=16
\ssbsectno=0
\def\section#1{\advance \sectno by 1 \sbsectno = 0
\sbsectno=0\ssbsectno=0\vskip 24pt
{\goodbreak
\noindent {\number\sectno.\ {#1}}}
\nobreak\vskip 12pt }
\def\subsection#1{\advance \sbsectno by 1
\ssbsectno=0
{\ifnum\count12=0\nobreak \else \medbreak \fi
\vskip 12 pt
\noindent {\number\sectno.\number\sbsectno. {\bf #1}}}
\nobreak\vskip 12pt }
\def\subsubsection#1{\advance \ssbsectno by 1
\medbreak\vskip 12pt{
\noindent{\number\sectno.\number\sbsectno.\number\ssbsectno.
\it #1.\ }}}
%
%
\countdef\figno=13
\figno = 0
\def\fig#1#2{\advance\figno by 1 \tolerance=10000
\setbox0=\hbox{\rm Fig. 00\ }
{\topinsert
\vskip #1 \smallskip\hangindent\wd0%
\noindent Fig.~\the\figno .\ #2
\vskip 12pt\endinsert}}
%
%
%

%
%
\def\REFERENCES
     {\countdef\refno=14
     \refno = 0\vskip 24pt
     {\centerline{\bf REFERENCES}\par\nobreak}
      \parindent=0pt
      \vskip 12pt\frenchspacing}
\def\ref#1#2#3
   {\advance\refno by 1
    \item{\rm[\number\refno] }\rm #1\ {\sl #2}
    \rm #3\par}
%
%
\def\frac#1#2{{#1 \over #2}}
%
%

%
%

%
%

%
%
\def\i{\ifmmode{\rm i}\else\char"10\fi}
%
%

%
\newfam\bofam
\font\tenbo=cmmib10   \textfont\bofam=\tenbo

\mathchardef\Omega="710A
\mathchardef\alpha="710B
\mathchardef\beta="710C
\mathchardef\gamma="710D
\mathchardef\delta="710E
\mathchardef\epsilon="710F
\mathchardef\rho="711A
\mathchardef\sigma="711B
\mathchardef\tau="711C
\mathchardef\upsilon="711D
\mathchardef\phi="711E
\mathchardef\chi="711F
\mathchardef\Gamma="7100
\mathchardef\Delta="7101
\mathchardef\Theta="7102
\mathchardef\Lambda="7103
\mathchardef\Xi="7104
\mathchardef\Pi="7105
\mathchardef\Sigma="7106
\mathchardef\Upsilon="7107
\mathchardef\Phi="7108
\mathchardef\Psi="7109
\mathchardef\zita="7110
\mathchardef\eta="7111
\mathchardef\theta="7112
\mathchardef\iota="7113
\mathchardef\kappa="7114
\mathchardef\lambda="7115
\mathchardef\mu="7116
\mathchardef\nu="7117
\mathchardef\xi="7118
\mathchardef\pi="7119
\mathchardef\psi="7120
\mathchardef\omega="7121
\mathchardef\varepsilon="7122
\mathchardef\vartheta="7123
\mathchardef\varpi="7124
\mathchardef\varrho="7125
\mathchardef\varsigma="7126
\mathchardef\varphi="7127

\baselineskip=18pt
\hsize=16truecm
\vsize=22truecm
\centerline{\bf ANALYTIC PROPERTIES OF}
\centerline{\bf THE RETURN MAPPING}
\centerline{\bf OF LI\'ENARD EQUATIONS.}
\bigskip\bigskip\bigskip
\centerline{\it J.P. Fran\c{c}oise}
\bigskip
\centerline{Universit\'e de Paris VI,} \centerline{Laboratoire
"GSIB",} \centerline{UFR de Math\'ematiques,}
 \centerline{175 rue de Chevaleret} \centerline{75013 Paris, France.}
\centerline{e-mail: jpf@ccr.jussieu.fr}
\bigskip
 \bigskip \bigskip\bigskip
\bigskip\bigskip

{\bf INTRODUCTION}
{\vskip 12pt}

         Bautin approach to the bifurcation theory of
limit cycles has been recently generalized in the framework of
complex analysis ([8][9]).
There are now more cases where the Bautin ideal is known.
A systematic study of the Poincar\'e center focus problem via
Abel equations entailed several new examples of Bautin
ideals ([2],[3]).
        This article deals with Li\'enard equations which have
been used in many applications (cf. [10]).
 Li\'enard equations play certainly a key role in Hilbert's
$16^{th}$ problem, as suggested for instance in ([19],[20]),
 because of the topological simplicity of the
return mapping . Limit
cycles encircle the origin and are necessarily contained in
the domain of existence of the return mapping. This domain
of existence may of course not be equal to the domain of convergence of the
analytic series which defines the return mapping in a neighborhood
of the origin. Nevertheless it is interesting to produce an
estimate of the size of this domain of convergence.
        Many contributions have been previously done to Li\'enard
equations and particularly by N.G. Lloyd and his co-workers ([1],[5],[14],[15]
,[16], [17]). But the
approach via the Lyapunoff series which is used in all these
references does not entail information on the domain of
convergence of the first return mapping itself.
        We develop here a direct computation of the return mapping
(indeed of its converse) which yields an estimate of the domain of convergence by a
majorizing series techniques and a recurrency relation for the coefficients
which entail the Bautin ideal.
        Application of
the complex analysis methods of ([8],[9]) is then straightforward. In Li\'enard case,
the use of the Hironaka  polynomial division theorem
is replaced by a very simple and explicit argument.
This displays a
bound to the number of complex limit cycles in a fixed
neighborhood of the origin.
        The maximal number of complex limit cycles on this disk equals $2n$.
Examples produced with perturbation theory (first-order)
of Hamiltonians by A. Lins Neto, W. de Melo and C.C. Pugh ([12]) show that
this bound is optimal
(recall that the complex bound is at least twice the real bound).
        To end this introduction we mention other recent references to
generalized Li\'enard equations([6],[7],[11],[13],[18])
to which eventually the same type of analysis should apply.
        The author thanks M. Caubergh for a careful reading
of the manuscript. This article was completed during a visit to
CRM of Montr\'eal. The author is grateful to L. Belair, C.Rousseau
and D. Schlomiuck for their invitation and to V. Kaloshin and S.
Yakovenko for useful comments related to the presentation of the
results.

{\vskip 12pt}

{\bf I-COMPUTATION OF THE BAUTIN IDEAL}
{\vskip 12pt}
        The perturbation of the linear part into a focus adds one
real limit cycle to the number obtained when the linear part is a
center. We consider only the case where the linear part is a
center to simplify the computations.
        This paragraph is concerned with the polynomial vector fields $X$ of
following type

$$X=x{\partial}/{\partial y}-y{\partial}/{\partial x}+
 \sum_{i=1}^{d}[{\lambda}_{i} x^{i}]y{\partial}/{\partial y}. \eqno(1.1)$$
and the associated flow, solution of the system:

$$\dot{x}= y , \eqno(1.2a)$$
$$\dot{y}= -x +  \sum_{i=1}^{d}[{\lambda}_{i} x^{i}]y = -x+p(x)y. \eqno(1.2b)$$

        The equations of the flow yield a second order
differential equation classically named the Li\'enard equation. Note
that this differential system is sometimes written:

$$\dot{x}=Y-P(x), \eqno(1.2c)$$
$$\dot{Y}=-x, \eqno(1.2d)$$

which is easily changed into the preceding ones with $y=Y-P(x)$, if $P'(x)=-p(x)$.

        Write (1.2) in polar coordinates $(r,\theta)$:

 $$ x = r{\rm cos}{\theta},  y = r{\rm sin}{\theta}. \eqno(1.3)$$

        This displays:

 $$ 2r{\dot r} = 2(x{\dot x}+y{\dot y}), r{\dot r}= r^{2}p(r{\rm cos}\theta)
 {\rm sin}^{2}\theta \eqno(1.4)$$

 $$ {\dot \theta} = (x{\dot y}-y{\dot x})/(x^{2}+y^{2}) = -1+{\rm sin}{\theta}{\rm cos}
 {\theta}p(r{\rm cos}\theta). \eqno(1.5)$$

        This yields:

$$dr/d{\theta}=rp(r{\rm cos}{\theta}){\rm sin}^{2}{\theta}/[
-1+{\rm sin}{\theta}{\rm cos}{\theta}p(r{\rm cos}{\theta})]. \eqno(1.6)$$

  Bautin's approach is based on the study of
    solutions of (1.6) $r = r(\theta)$ so that $r(0) = r_{0}$, given as
   an expansion:

 $$r_{0} = r + v_{2}(\theta)r^{2}+ ... +  v_{k}(\theta)r^{k}+ ...\eqno(1.7)$$

 Comparison between (1.6) and (1.7) yields:

$$\sum_{k \geq 1}v_{k}'(\theta)r^{k}
[-1+{\rm sin}{\theta}{\rm cos}{\theta}p(r{\rm cos}{\theta})]
$$
$$+\sum_{k \geq 1}[kv_{k}(\theta)p(r{\rm
cos}{\theta}){\rm sin}^{2}{\theta}]r^{k}=0. \eqno(1.8)$$

This displays the following recurrency relation on the
coefficients $v_{k}(\theta)$:

$$v_{k}'(\theta)=\sum_{l=1}^{d} {\lambda}_{l}({\rm cos}{\theta})^{l}{\rm
sin}{\theta}[{\rm cos}{\theta}v'_{k-l}(\theta)+(k-l){\rm
sin}{\theta}v_{k-l}(\theta)]. \eqno(1.9)$$

        Let $I$ be an ideal of $R[{\lambda}_{1},...,{\lambda}_{d}]$. It is
convenient to denote $v_{k}(\theta)\in I$ to mean that for all
values of $\theta$, $v_{k}(\theta)$ is a polynomial in the
parameters $({\lambda}_{1},...,{\lambda}_{d})$ which belongs to the ideal $I$. In
the following $c_{k}$ denotes a sequence of
non-zero numbers (independent of the parameters) which will be defined
inductively.
        Choose the initial conditions for the recurrency relation (1.9)
as:

$$ v_{1}(\theta)=1. \eqno(1.10)$$

        Then (1.9) yields:
$$v'_{2}(\theta)={\lambda}_{1}{\rm cos}{\theta}{\rm sin}^{2}{\theta},
 {\rm thus}: v_{2}(\theta)={\lambda}_{1}w_{2}^{1}({\rm sin}{\theta}),
 w_{2}^{1}(0). \eqno(1.11a)$$
        This entails:
 $$v_{2}(2\pi)=0. \eqno(1.11b)$$
        The first coefficient which really matters for the Bautin
ideal is the next one. The recurrency relation (1.9) displays:

$$v'_{3}(\theta)={\lambda}_{1}{\rm cos}{\theta}{\sin}{\theta}[{\rm cos}{\theta}v'_{2}
+2{\rm sin}{\theta}v_{2}]+{\lambda}_{2}{\rm cos}^{2}{\theta}{\rm
sin}^{2}{\theta}=0. \eqno(1.12)$$

        This yields:

$$v_{3}(\theta)= f_{3}({\theta})+{\lambda}_{1}w_{3}^{1}({\rm sin}{\theta}),\eqno(1.13a)$$
with:
$$w_{3}^{1}(0)=0, \eqno(1.13b)$$
and
$$f_{3}({\theta})\in ({\lambda}_{2}), \eqno(1.13c)$$

and
$$v_{3}(2\pi)=c_{3}{\lambda}_{2}, c_{3}={\int}_{0}^{2\pi}{\rm
cos}^{2}{\theta}{\rm
sin}^{2}{\theta}d{\theta}\neq 0. \eqno(1.14)$$

        We prove now by induction the following:

{\vskip 6pt}
{\bf Lemma I.1}
{\vskip 6pt}
        Let $k_{0}$ be the maximal integer so that $2k_{0}\leq k-1$ and $k_{1}$ be the
maximal integer so that $2k_{1}+2\leq k$.
        The coefficient $v_{k}(\theta)$ displays the decomposition:

$$v_{k}(\theta)=
f_{k}(\theta)+\sum_{j=0}^{k_{1}}{\lambda}_{2j+1}w_{k}^{2j+1}({\rm
sin}{\theta}, {\lambda}), \eqno(1.17a)$$

where:

$$f_{k}(\theta)\in ({\lambda}_{2},...,{\lambda}_{2k_{0}}), \eqno(1.17b)$$

$$w_{k}^{2j+1}(0, \lambda)=0. \eqno(1.17c)$$

{\vskip 6pt}
{\bf Proof}
{\vskip 6pt}
        This is certainly true for $v_{2}(\theta)$ and
$v_{3}(\theta)$. Assume this is so inductively. The recurrency
relation (1.9) shows that the term $w_{k-l}^{2j+1}({\rm
sin}{\theta})$ contributes either (for $l$ even) to an element of
the ideal $({\lambda}_{2},..., {\lambda}_{2k_{0}})$ or (for $l=2h+1$ odd) produces:

$${\lambda}_{2h+1}({\rm cos}{\theta})^{2h+1}{\rm sin}{\theta}[{\rm
cos}^{2}{\theta}w_{k-l}^{2j+1}({\rm sin}{\theta})'+{\rm
sin}{\theta}w_{k-l}^{2j+1}({\rm sin}{\theta})], \eqno(1.18)$$

which (once integrated against $\theta$ contributes to
${\lambda}_{2h+1}w_{k}^{2h+1}({\rm sin}{\theta}))$.

{\vskip 6pt}
        An immediate consequence of the proposition is the:
{\vskip 6pt}
{\bf Lemma I.2}
{\vskip 6pt}
        The coefficients $v_{k}(2\pi)$ belong to the ideal
${\lambda}_{2},...,{\lambda}_{2k_{0}}$.
{\vskip 6pt}
        This last result can be improved for the coefficients
$v_{k}(2\pi)$ of odd order $k$. Denote $n=[d/2]$ (integer part of $d/2$).
{\vskip 6pt}
{\bf Lemma I.3}
{\vskip 6pt}
        For all odd values of $k=2k_{0}+1$ $(k_{0}=1, ... , n)$,
 the coefficient

$v_{2k_{0}+1}(2\pi)$ is such that:

$$v_{2k_{0}+1}(2\pi){\in}({\lambda}_{2},...,{\lambda}_{2k_{0}-2})+
c_{2k_{1}+1}{\lambda}
_{2k_{0}},
\eqno(1.19)$$

with:

$$c_{2k_{0}+1}=\int_{0}^{2\pi}({\rm cos}{\theta})^{2k_{0}}({\rm
sin}{\theta})^{2}d{\theta}\neq 0. \eqno(1.20)$$
{\vskip 6pt}
{\bf Proof}
{\vskip 6pt}
        In the recurrency relation (1.9), the only term which
contributes to ${\lambda}_{2k_{0}}$ is:

$${\lambda}_{2k_{0}}({\rm cos}{\theta})^{2k_{0}}{\rm
sin}{\theta}[{\rm
sin}{\theta}],\eqno(1.21)$$

which yields (1.19).

        The recurrency relation above yields finally:

{\vskip 6pt}
{\bf Theorem 1.3}
{\vskip 6pt}
        The Bautin ideal of the Li\'enard vector fields (1.1),
defined as the ideal generated by the coefficients of the
return mapping (or of its inverse) $v_{k}(2\pi)$ is equal to the
ideal $({\lambda}_{2},..., {\lambda}_{2n})$.
 The Bautin index $B$ defined as the first
integer $k$ so that the increasing sequence of ideals generated by
the $k$ first coefficients $v_{k}(2\pi)$ becomes stationary is
equal to $2n+1$.
{\vskip 12pt}
        The recurrency relation entails as well the following
result:
{\vskip 12pt}
{\bf Theorem 1.4}
{\vskip 12pt}
        For all values of $\lambda$, there is a neighborhood of
the origin on which the number of real limit cycles is less than
$n-1$ ($n$ if we consider the perturbation of a focus).
{\vskip 12pt}
{\bf Proof}
{\vskip 12pt}
        This is a consequence of the classical Bautin's argument.
Collecting the terms of the first return mapping, we write the
equation for the real limit cycles as:

$$v_{3}(2\pi)r^{3}(1+...)+v_{5}(2\pi)r^{5}(1+...)+...+v_{2n+1}(2\pi)r^{2n+1}(1+...)
=0.
\eqno(1.22)$$
        Successive applications of Rolle's lemma show that the number of real
positive zeroes of (1.22)is
less than $n-1$. This shows that the maximal number of limit cycles
 which can bifurcate
when the linear part is of focus type is $n$ (because it adds up one term in (1.22).
Note that this bound is in agreement with the
bound foreseen in A. Lins neto-W. de Melo-C.C. Pugh conjecture.
(cf [12]). Of course this computation does not entail any control
of the size of the domain on which the number of limit cycles is
less than $n$ in terms of the coefficients of the perturbation.
This is the reason to develop further analysis based on the
complexification of the return mapping.
{\vskip 12pt}
{\bf II- Estimates of the radius of convergence of the first return
mapping}
{\vskip 6pt}

 Let $ f_{\lambda}(x)= \sum a_{k}(\lambda)x^{k}$ be an analytic series in $x$ with
polynomial coefficients in the parameters $\lambda= (\lambda_{1},...,\lambda_{d})$.
 Denote
$\mid a_{k}\mid$ (norm of the polynomial $a_{k}$) as the sum of
the absolute
value of the coefficients and $\mid\lambda\mid=\mid\lambda_{1}\mid+...
+\mid\lambda_{d}\mid$. Recall now the following:
\vskip 12pt
{\bf Definition II.1}
\vskip 12pt

        The series $f_{\lambda}$ is called an $A_{0}$-series if the following
 two conditions are satisfied:

 There are positive constants $K_{1}, K_{2},K_{3},K_{4}$ such that:

 $$deg(a_{k})\leq K_{1}k+K_{2},\eqno(2.2a)$$

 $$\mid a_{k}\mid \leq K_{3}K_{4}^{k}.\eqno(2.2b)$$
 \vskip 12pt

 $A_{0}$-series form a subring of the ring of formal power series in $x$ with polynomial
coefficients in $\lambda$. All the usual analytic operations, like substitution
to a given analytic
function, composition, inversion,... transform $A_{0}$-series into themselves.
 $A_{0}$-series have been precisely introduced (in the
subject) by M. Briskin and Y. Yomdin ([4]).

\vskip 12pt
        In the following, we also denote  $f_{\lambda}$ the complex analytic
function defined for all
$\lambda \in C^{D}$ on a disc $D(0,R)$ by the $A_{0}$-series.

        In Li\'enard case, the following holds:
{\vskip 12pt}
{\bf Proposition II.2}
{\vskip 12pt}
        The series (1.7) is for all values of $\theta$ a
$A_{0}$-series with:
$$K_{1}=1, K_{2}=-1, K_{3}=\pi/2,K_{4}=2.\eqno(2.4)$$
{\vskip 12pt}
{\bf Proof.}
{\vskip 12pt}
        The inverse series of (1.7), which is indeed the
Bautin series:

$$r=r_{0}+w_{2}(\theta)r_{0}^{2}+...+w_{k}(\theta)r_{0}^{k}+...\eqno(2.5)$$

solves the differential equation (1.6).
        We use Siegel majorizing series techniques. A rough
majorizing series (both in ${\lambda}$ and $r_{0}$)
of (1.6) is provided by the solution $M$ of the differential
equation:

$$dM/d{\theta}=[M{\mid\lambda\mid}M/(1-M)]/[1-{\mid\lambda\mid}M/(1-M)],\eqno(2.6)$$

which displays the expansion:

$$M=r_{0}+M_{2}(\theta)r_{0}^{2}+...+M_{k}(\theta)r_{0}^{k}+...\eqno(2.7)$$

        The coefficients $M_{k}(\theta)$ are inductively defined
by a recurrency relation which is of following type:

$$M'_{k}(\theta)= S_{k}[M_{2}(\theta),...,M_{k-1}(\theta)], k \geq
2,\eqno(2.8)$$

where $S_{k}$ is a polynomial with positive coefficients.

        Denote
$$M_{k}={\rm Max}[M_{k}(\theta), {\theta}\in [0, 2\pi]].
\eqno(2.9)$$

The equation (2.8) yields the following inequality:

$$M_{k} \leq (2\pi)S_{k}[M_{2},...,M_{k-1}], k\geq 2. \eqno(2.10)$$

        This displays a new majorizing series for (2.5) $W(r_{0})$
solution of the algebraic equation:

$$W(r_{0})-r_{0}=(2\pi){\mid\lambda\mid}W(r_{0})^{2}/[1-(1+{\mid\lambda\mid})W(r_{0})].
\eqno(2.11a)$$

        The algebraic equation (2.11a) has a unique analytic
solution $W(r_{0})$ which is tangent to $r_{0}$ for small values
of $r_{0}$. This equation yields $r_{0}$ in terms of $W(r_{0})$:

$$r_{0}=W(r_{0})-(2\pi){\mid\lambda\mid}W(r_{0})^{2}/[1-(1+{\mid\lambda\mid})
W(r_{0})].
\eqno(2.11b)$$

        Elementary considerations on majorizing series show that the converse
of a majorizing series provides a majorizing series for the converse.

        This immediately provides an estimate of the radius of convergence
 of (1.7):

$$R(\lambda)= 1/[1+\mid\lambda\mid], \eqno(2.12)$$

and the proof of the proposition II-2.

        The estimate of the radius of convergence can be improved
(for some values of the parameters) with an elementary scaling
argument in the case of Li\'enard's equations. The proof given
above works actually for an arbitrary polynomial perturbation of
the rotation flow (with some changes of notations).
{\vskip 12pt}
{\bf Theorem II.3}
{\vskip 12pt}
        The return mapping of the vector field $X$ converges at
least on the disc $D(0, R(\rho))$ of radius $R(\rho)=\rho/2$ where
$\rho$ is the unique positive real number such that:

$$\rho^{d}\mid\lambda_{d}\mid+...+\rho\mid\lambda_{1}\mid=1.
\eqno(2.13)$$

        This last estimation improves (2.12) for small
$\mid\lambda\mid$ but not for large ones.
{\vskip 6pt}
Proof:
{\vskip 6pt}
        Change of coordinates $(x,y)$ into $({\rho}x,{\rho}y)$
transforms

$$\dot{x}=y$$
$$\dot{y}=-x+p(x)y,\eqno(2.14a)$$

into

$$\dot{x}=y$$
$$\dot{y}=-x+\sum_{i=1}^{d}(\lambda_{i}\rho^{i})x^{i}y.
\eqno(2.14b)$$

        This means that any result obtained on the disc $D(0,R)$
for the equation (2.14a) is valid on the disc
$D(0,{\rho}R)$ for the new equation $\lambda_{i}\mapsto
        \lambda'_{i}=\rho^{i}\lambda_{i}$.
{\vskip 12pt}
{\bf III-Estimates of the number of complex limit cycles on a
fixed neighborhood of the origin.}
{\vskip 12pt}

        The generalization of Bautin's approach to complex
limit cycles has been done in ([8],[9]). It yields quite explicit
results in the case of homogeneous perturbations of the linear
part. The same type of results can be displayed in the case of
Li\'enard equations and this is done in this paragraph.

        It is quite interesting to note now that in comparison to
the "polynomial Hironaka division theorem" that we need to
use in the general situation, we have to check a very
easy proposition.

        In Li\'enard case, the Hironaka basis is obviously
${\lambda}_{2},...{\lambda}_{2n}$ and it yields the
following:

{\vskip 12pt}
{\bf Proposition III.1}
{\vskip 12pt}
        Let $I$ be the Bautin ideal. Let $f({\lambda})$ be an
element of $I$ of degree $k$, there is a decomposition:
$$f(\lambda) = \sum_{i=1}^{n} \phi_{i}(\lambda)\lambda_{2i}, \eqno(3.1a)$$

with

$$deg(\phi_{i}) \leq deg(f)-1 =k-1, \eqno(3.1b) $$

and

$$\mid \phi_{i}\mid \leq \mid f\mid. \eqno(3.1c) $$
{\vskip 12pt}
{\bf Proof.}
{\vskip 12pt}
        Consider $f(\lambda)$, collect in front of $\lambda_{2}$
all the monomials which contain $\lambda_{2}$. The difference
still belongs to the ideal $I$. Then repeat the process with
${\lambda}_{4}$,...The two majorations are obvious.

        Following the techniques of ([8],[9]), we prove now that the
first return mapping belongs to a Bernstein class. It is
convenient at this point to change notations and write:

$$f(\lambda, r)=r+v_{2}(2\pi)r^{2}+...
+v_{k}(2\pi)r^{k}+...,\eqno(3.2)$$
for the (converse of the) first return mapping and
 $f_{k}(\lambda)=v_{k}(2\pi)$ for its
coefficients. Indeed, we consider the analytic extension of the
series to the complex domain and denote now $x$ as the complex
variable in place of the real variable $r$.

$$f(\lambda, x)=x+
f_{2}(\lambda)x^{2}+...+f_{k}(\lambda)x^{k}+...\eqno(3.3)$$

{\vskip 12pt}
{\bf Proposition III-2}
{\vskip 12pt}
        The analytic series (3.3) belongs to the Bernstein class
$B_{(B, R, c)}^{2}$ with:

$$R=K_{4}(1+\mid\lambda\mid),\eqno(3.4a)$$
$$c=nK_{3}^{2}K_{4}^{2n-1}/c_{2n+1}^{n}.
\eqno(3.4b)$$

{\vskip 12pt}
Proof.
{\vskip 12pt}
 Write first the condition for $f_{\lambda}(x)$ to be an
$A_{0}$-series as follows:
$$deg[f_{k}]\leq k-1,$$
$$\mid f_{k}\mid \leq K_{3}K_{4}^{k}.$$

Write next the decomposition:

$$f_{k}(\lambda)=\sum_{i=1}^{n}[\phi_{k, i}(\lambda)\lambda_{2i}]. \eqno(3.5)$$

This yields:

$$\mid \phi_{k, i}(\lambda)\mid \leq \mid
f_{k}\mid(1+\mid\lambda\mid)^{k-2}\leq K_{3}K_{4}^{k-1}(1+\mid\lambda\mid)]^{k-2}.
\eqno(3.6)$$

This entails:

$$\mid f_{k}(\lambda \leq K_{3}K_{4}^{k-1}(1+\mid\lambda\mid)]^{k-2}({\rm
Max}_{i=1,...,n}\mid\lambda_{2i}\mid). \eqno(3.7)$$

A more careful analysis of the recurrency relation (cf. 1.19, 1.20)
shows that equation (3.5) entails a $nxn$ matrix relation
between the vectors $f_{2j+1}, j=1,...,n$ and $\lambda_{2i},
i=1,...,n$ of the form:

$$f_{2j+1}=\sum_{i=1}^{n}[C_{ji}+\Delta_{ji}]\lambda_{2i},
\eqno(3.8)$$

where the matrices $C$ and $\Delta$ are respectively diagonal, with
non-zero coefficients $c_{2i+1}$ defined in (1.20), and nilpotent (upper-triangular).
Inverting the matrix relation (3.8) yields:

$$\mid\lambda_{2i}\mid \leq \mid(C+\Delta)^{-1}\mid{\rm
Max}_{i=1,...,n}\mid f_{2i+1}(\lambda)\mid. \eqno(3.9)$$

This is completed with the inequalities:

$$\mid (C+\Delta)^{-1}\mid \leq
\mid C^{-1}\mid\mid(1+C^{-1}\Delta)^{-1}\mid, \eqno(3.10a)$$

$$\leq \mid C^{-1}\mid[1+\mid C^{-1}\Delta\mid+...+(\mid C^{-1}\Delta\mid)^{n-1}].
\eqno(3.10b)$$

The coefficients of the diagonal matrix $C$ entail (cf. 1.20):

$$\mid C^{-1}\mid=1/c_{2n+1}, \eqno(3.11)$$

This yields:

$$\mid (C+\Delta)^{-1}\mid \leq
(1/c_{2n+1})^{n}nK_{3}K_{4}^{2n}(1+\mid\lambda\mid)]^{2n+1}. \eqno(3.12)$$

Equations (3.12), (3.9), (3.7) entail now:

$$\mid f_{k}(\lambda)\mid R^{k}\leq (1/c_{2n+1})^{n}nK_{3}^{2}
[K_{4}(1+\mid\lambda\mid)]^{2n-1}
{\rm Max}_{k=2,..., B}(\mid f_{k}(\lambda)\mid), \eqno(3.13a)$$

$$\mid f_{k}(\lambda)\mid R^{k}\leq
(1/c_{2n+1})^{n}nK_{3}^{2}[K_{4}(1+\mid\lambda\mid)]^{4n}
{\rm Max}_{j=2,..., B}(\mid f_{j}(\lambda)\mid)R^{j}. \eqno(3.13b)$$

This means that the analytic series (3.3) belongs to the Bernstein class
$B_{(B, R, c)}^{2}$ (cf. [8], [9])with:
$$c=nK_{3}^{2}[K_{4}(1+\mid\lambda\mid)]^{4n}/c_{2n+1}^{n}.$$
\vskip 12pt
        The proposition III.2 and the results of ([8], [9])
now imply the following:
\vskip 12pt
{\bf Theorem III.3}
\vskip 12pt
        The number of zeros of $f_{\lambda}(x)$ in the disc $D(0, R')$
is less than $B-1=2n$ with
$$R'=c_{2n+1}^{n}/[2^{6n}nK_{3}^{2}[K_{4}(1+\mid\lambda\mid)]^{4n+1},\eqno(3.14a)$$

$$R'=c_{2n+1}^{n}/[2^{10n-1}n{\pi}^{2}(1+\mid\lambda\mid)^{4n+1}.\eqno(3.14a)$$

        For small values of $\mid\lambda\mid$, this estimate can
be improved with the same scaling argument as used in the previous
paragraph as follows:

{\vskip 12pt}
{\bf Theorem III.4}
{\vskip 12pt}
        The vector field $X$ has less than $2n$ complex limit
cycles on the disc $D(0, R'(\rho))$ of radius
$$R'(\rho)={\rho}c^{n}_{2n+1}/[{\pi}^{2}n2^{14n}], \eqno(3.15a)$$

where
$\rho$ is the unique positive real number such that:

$$\rho^{d}\mid\lambda_{d}\mid+...+\rho\mid\lambda_{1}\mid=1.
\eqno(3.15b)$$

\vfill\eject
\REFERENCES

\ref{T.R. Blows, N. G. Lloyd:}{}{The number of small-amplitude
limit cycles of Li\'enard equations. Math. Proc. Cambridge Philos.
Soc. {\bf 95}, 359-366 (1984).}

\ref{M. Briskin, J.-P. Francoise, Y. Yomdin:}{}{The Bautin ideal
of the Abel Equation. Nonlinearity {\bf 11}, 431-443 (1998).}

\ref{M. Briskin, J.-P. Francoise, Y. Yomdin:}{}{Center conditions,
composition of polynomials and moments on algebraic curves.
Ergod. Th.  Dynam. Sys. {\bf 19}, 1201-1220 (1999).}

\ref{M. Briskin, Y. Yomdin:}{}{Algebraic Families of Analytic
Functions. J. Diff. Equations. {\bf 136} (2), 248-267 (1997).}

\ref{C.J. Christopher, N.G. Lloyd:}{}{Small-amplitude limit cycles
in Li\'enard systems. Non-linear Diff. Equations and Applications.
{\bf 3}, 183-190, (1996).}

\ref{G. Dangelmayr, J. Guckenheimer:}{}{On a four parameter family
of planar vector fields. Arch. Ration. Mech. {\bf 97}, 321-352,
(1987).}

\ref{F. Dumortier, C. Rousseau:}{}{Cubic Li\'enard equations with
linear damping. Nonlinearity. {\bf 9}, 1489-1500, (1990).}

\ref{J.-P. Francoise, Y. Yomdin:}{}{ Bernstein inequality and applications
to analytic geometry and differential equations. J. Funct. Analysis {\bf 146},
 185-205 (1997).}

\ref{J.-P. Francoise, Y. Yomdin:}{}{Projection of analytic sets and Bernstein
inequalities. Singularities Symposium-Lojasiewicz 70, Edts B. Jacubczyk,
W. Pawlucki, Y. Stasica, Banach Center Publications, Warszawa {\bf 44},
103-108 (1998).}

\ref{M. Hirsh, S. Smale:}{}{{\it Differential Equations, Dynamical
Systems and Linear Algebra}. Academic Press, New York.}

\ref{A. Gasull, J. Torregrosa:}{}{Small-amplitude limit cycles in
Li\'enard systems via multiplicity. J. Diff. Equations}

\ref{A. Lins neto, W. de Melo, C.C. Pugh:}{}{On Li\'enard
equations. Geometry and Topology (Rio de Janeiro, 1976). Lecture
Notes on Mathematics, {\bf 597}, 335-357, (1977).}

\ref{A.I. Khibnik, B. Krauskopf, C. Rousseau:}{}{Global study of a
family of cubic Li\'enard equations. Nonlinearity {\bf 11},
1505-1519 (1998).}

\ref{N.G.Lloyd:}{}{Li\'enard systems with several limit cycles.
Math. Proc. Cambridge Philos. Soc. {\bf 102}, 565-572, (1987).}

\ref{N.G. Lloyd, S. Lynch:}{}{Small-amplitude limit cycles of
certain Li\'enard systems. Proc. R. Soc. Lond. A {\bf 418},
199-208 (1988).}

\ref{S. Lynch:}{}{Small-amplitude limit cycles of Li\'enard
systems. Calcolo {\bf 27}, 1-32, (1990).}

\ref{S. Lynch:}{}{Li\'enard systems and the second part of
Hilbert's $16^{th}$ problem. Non linear Anal. {\bf 30}, 1395-1403,
(1997).}

\ref{S.Malo:}{}{Rigorous computer verification of planar vector
field structure. {\it PhD Thesis} Cornell University (1994).}

\ref{S. Smale:}{}{Dynamics retrospective: great problems, attempts
that failed. Physica D, {\bf 51}, 267-273, (1991).}

\ref{S. Smale:}{}{Mathematical Problems for the Next Century. The
Mathematical Intelligencer, {\bf 20}, 7-15, (1998).}

\end